\input amstex
\documentstyle{amsppt}
\magnification 1200
\vcorrection{-9mm}
\NoBlackBoxes
\input epsf

\let\wt\widetilde
\def\Z{\Bbb Z}
\def\R{\Bbb R}
\def\C{\Bbb C}
\def\P{\Bbb P}    \let\pp\P
\def\sph{\Bbb S}
\def\RP{\Bbb{RP}} \let\rp\RP
\def\CP{\Bbb{CP}}
\def\conj{\operatorname{conj}}

\def\hyp{\operatorname{hyp}}
\def\ord{\operatorname{ord}}
\def\lk{\operatorname{lk}}

\def\supp{\operatorname{supp}}
\def\eps{\varepsilon}
\def\Tproj{T_{\operatorname{proj}}}
\def\cn{\operatorname{cr}}
\def\ps{\operatorname{ps}}

\def\rl{l}

\topmatter
\title    Maximally writhed real algebraic links
\endtitle

\author   Grigory Mikhalkin and Stepan Orevkov
\endauthor

\abstract
             Oleg Viro introduced an invariant of rigid isotopy for real
             algebraic knots and links in $\Bbb{RP}^3$ which can be viewed
             as a first order Vassiliev invariant.
             In this paper we classify real algebraic links of degree $d$
             with the maximal value of this invariant in its two versions:
             $w$ and $w_\lambda$.
\endabstract

\thanks
Research is supported in part by  the SNSF-grants
159240, 159581 and the NCCR SwissMAP project  (G.M.)
\endthanks

\address        Universit\'e de Gen\'eve, Section de Math\'ematiques,
                Battelle Villa, 1227 Carouge, Suisse.
\endaddress
\email          grigory.mikhalkin\@unige.ch
\endemail

\address
                Steklov Mathematical Institute, Gubkina 8, 119991, Moscow, Russia;\smallskip
                IMT, Universit\'e Paul Sabatier,
                118 route de Narbonne, 31062, Toulouse, France.
\endaddress
\email          orevkov\@math.ups-tlse.fr
\endemail

\endtopmatter

\def\refB    {1}  
\def\refD    {2}  
\def\refEM   {3}  
\def\refGM   {4}  
\def\refGH   {5}  
\def\refGroH {6}  
\def\refMa   {7}  
\def\refMO   {8}  
\def\refMOd  {9} 
\def\refMu   {10} 
\def\refN    {11} 
\def\refRo   {12} 
\def\refV    {13} 

\def\sectHopf        {2}

\def\sectProofMainTh {3}
\def\sectTors        {\sectProofMainTh.1}
\def\sectHyperbPla   {\sectProofMainTh.2}
\def\sectHyperbSpa   {\sectProofMainTh.3}  \let\sectHypSpa\sectHyperbSpa
\def\sectT           {\sectProofMainTh.4}
\def\sectEOP         {\sectProofMainTh.5}
\def\sectPropU       {\sectProofMainTh.6}

\def\sectPermut      {4}
\def\sectHabTpq      {\sectPermut.1} \let \sectTpqHab \sectHabTpq
\def\sectPS          {\sectPermut.2}
\def\sectDiagr       {\sectPermut.3}

\def\sectPermutLast  {\sectPermut.6}

\def\sectExist       {5}

\def\thMain       {1}
\def\propPS       {1.1}
\def\propMurasugi {1.2}
\def\conjMurasugi {1.3}
\def\conjRigid    {1.4}
\def\thAnyGenus   {2}
\def\corAnyGenus  {1.6}
\def\remWO        {1.5}
\def\thExist      {3}

\def\propHopf       {\sectHopf.1}

\def\lemTors        {\sectProofMainTh.1}
\def\lemHypC        {\sectProofMainTh.2}
\def\lemHypTypeI    {\sectProofMainTh.3}
\def\lemHypKC       {\sectProofMainTh.4}
\def\lemHypK        {\sectProofMainTh.5}
\def\lemHypPosA     {\sectProofMainTh.6}
\def\lemHypPos      {\sectProofMainTh.7}

\def\lemHypII       {\sectProofMainTh.8}
\def\lemHypIII      {\sectProofMainTh.9}
\def\lemDisjoint    {\sectProofMainTh.10}
\def\lemTwoLines    {\sectProofMainTh.11}
\def\lemLK          {\sectProofMainTh.12}
\def\propU          {\sectProofMainTh.13}
\def\propFibr       {\sectProofMainTh.14}
\def\remCycloid     {\sectProofMainTh.15}

\def\propTpqHab     {\sectPermut.1}  \let \propHabTpq \propTpqHab
\def\propSect       {\sectPermut.2}
\def\propLift       {\sectPermut.3}  
\def\thPermut       {\sectPermut.4}  
\def\remPermut      {\sectPermut.5}
\def\remNonPermut   {\sectPermut.6}
\def\lemCoxeter     {\sectPermut.7}

\def\lemExistPrtrb  {\sectExist.1}
\def\remExist       {\sectExist.2}
\def\lemNonSpe      {\sectExist.3}
\def\lemExistOne    {\sectExist.4}
\def\lemExistTwo    {\sectExist.5}

\def\eqDefWL    {1}
\def\eqDefHyp   {2}
\def\eqSumA     {3}
\def\eqCRW      {4}
\def\eqT        {5}

\def\eqNonSpe   {6}

\def\figHopf    {1}
\def\figHypPos  {2}
\def\figCusp    {3}
\def\figChu     {4}
\def\figChuInv  {5}
\def\figTwoChu  {6}
\def\figSchPrf  {7}
\def\figR       {8}
\def\figCycloid {9}
\def\figT       {10}
\def\figSmoo    {11}
\def\figSmooG   {12}

\def\CondW  {(i)}
\def\CondWA {(i$'$)}
\def\CondH  {(ii)}
\def\CondPS {(iii${}^{\text{t}}$)}
\def\CondPSA{(iii${}^{\text{a}}$)}
\def\CondCr {(iv${}^{\text{t}}$)}   
\def\CondCrA{(iv${}^{\text{a}}$)}   \let \CondCRA \CondCrA
\def\CondT  {(v)}
\def\CondTA {(v$'$)}

\document

\head 1. Introduction
\endhead

\subhead 1.1 \endsubhead
A {\it real algebraic curve} in $\P^3$ is a (complex) one-dimensional subvariety $A$ in
$\P^3=\CP^3$ invariant under the involution of complex conjugation
$\conj:\P^3\to\P^3$, $(x_0:x_1:x_2:x_3)\mapsto(\bar x_0:\bar x_1:\bar x_2:\bar x_3)$.
The conj-invariance is equivalent to the fact that $A$ can be defined by a system of
homogeneous polynomial equations with real coefficients.
The degree of $A$ is defined as its homological degree, i.~e. the number $d$ such that
$[A] = d[\P^1]\in H_2(\P^3)\cong\Z$. A curve of degree $d$ intersects a generic complex plane
in $d$ points.

We denote the real locus (i.e., the set of real points) of $A$ by $\R A$.
We say that a real curve $A$ is {\it smooth} if it is a smooth complex submanifold of $\P^3$.
In this case, $\R A$ is a smooth real submanifold of $\RP^3$ and if each of its irreducible components
is non-empty,
we call it a {\it real algebraic link} or, more specifically,
a {\it real algebraic knot} in the case when $\R A$ is connected.

According to our definition, if $L=\R A$ is an algebraic link, then $A$ is uniquely
determined by $L$. In this case we call it the {\it complexification of $L$} and we write
$A=\C L$. We emphasize again that, by definition, the complexification of a real algebraic link
is a {\sl smooth} complex subvariety of $\CP^3$.
We say that a real algebraic link is {\it irreducible} if its complexification is connected.
All links considered in this paper are irreducible.

Two real algebraic links are called {\it rigidly isotopic} if their complexifications belong to the same
connected component of the space of real algebraic links of the same degree.
A rigid isotopy classification of real algebraic knots of genus $0$ in $\P^3$ is obtained in [\refB] up to degree 5
and in [\refMO] up to degree 6. Also we gave in [\refMO] a rigid isotopy classification for genus one knots
and links up to degree 6 (here we speak of the genus of the complex curve $A$ rather than the
minimal genus of a Seifert surface of $\R A$).

It happens that in all the above-mentioned cases, a rigid isotopy class is
completely determined by the usual (topological) isotopy class, the complex orientation (for genus one links),
and the invariant of rigid isotopy $w$ introduced by Viro [\refV] (called in [\refV] {\it encomplexed writhe\/}).
This invariant is defined as the sum of signs of crossings of a generic projection but the real
crossings with non-real branches are also counted with appropriate signs, while the crossings
between different link components are not counted; see details in [\refV]
(the definition of $w$ is also reproduced in [\refMO]).

Let $T(p,q)=\{(z,w)\mid z^p=w^q\}\cap \sph^3$, $p\ge q\ge 0$,
be the $(p,q)$-torus link in the 3-sphere $\sph^3\subset\C^2$ of radius $\sqrt2$.
If $p\equiv q\mod2$, we define the {\it projective torus link} 
$\Tproj(p,q) = T(p,q)/(-1) \subset \sph^3/(-1) = (\C^2\setminus 0)/\R^{\times}=\RP^3$
(here we refer to an identification of $\R^4$ with $\C^2$).
Note that the link $\Tproj(p,q)$ sits in the torus $\{|z|^2=|w|^2=1\}/(-1)\subset\sph^3/(-1)$
which corresponds to a hyperboloid under the identification of $\sph^3/(-1)$ with $\RP^3$.
The link $\Tproj(p,q)$ is an algebraic curve of bidegree $(a,b)=((p+q)/2,(p-q)/2)$ on this hyperboloid,
see [\refMO; \S3] and Section \sectHabTpq.
This is why in [\refMO] we called this link {\it hyperboloidal link} $h_{a,b}$.
The curve $\Tproj(p,q)\subset\P^3$ is a real algebraic curve of degree $p$.
If $\gcd(a,b)=1$ then it is irreducible, and, furthermore, rational.
If $b=1$ then this curve is smooth (otherwise it has a pair of complex conjugate singularities
of type $u^a=v^b$), thus $\Tproj(d,d-2)$ is an algebraic knot.

\smallskip
For a link $L$ in $\RP^3$ (a differentiable real one-dimensional submanifold),
we define the {\it crossing number} $\cn(L)$ in the same way as for links in $\sph^3$. It is
the minimal number of crossings of a generic plane projection where the minimum is taken over
the isotopy class of $L$. Similarly we define the {\it plane section number} $\ps(L)$
as the minimal number of intersections with a generic plane where again the minimum is taken 
over the isotopy class of $L$.
The plane section number of a link in $\RP^3$ can be equivalently defined as the
smallest number of immersed arcs of any diagram of the link (in the sense of [\refD]; see also
Section \sectDiagr).

\smallskip
Let $N_d = (d-1)(d-2)/2$
(this is the maximal number of nodes for an irreducible planar curve
of degree $d$ and thus also the maximal possible value of $w$ for
real algebraic knots of degree $d$).
It is not difficult to check that $w(\Tproj(d,d-2))=N_d$.
The main result is the following.

\proclaim{ Theorem \thMain }
Let $K$ be a real algebraic knot of genus $0$ and of degree $d\ge 3$.
Then the following conditions are equivalent:
\footnote{The superscripts ${}^{\text{a}}$ and ${}^{\text{t}}$ in ``\CondPSA", ``\CondPS", etc.
come from `algebraic' and `topological'.}
\roster
\item"\CondW"   $|w(K)|=N_d$;
\item"\CondH"   any real plane tangent to $K$ has only real intersections with $\C K$;
\item"\CondPSA" any generic real plane cuts $K$ at $d$ or $d-2$ real points;
\item"\CondPS"  $\ps(K)=d-2$;
\item"\CondCrA" any generic projection of $K$ has at least $N_d-1$ hyperbolic double points;
\item"\CondCr"  $\cn(K)=N_d-1$;
\item"\CondT"   $K$ or its mirror image is isotopic to $\Tproj(d-2,d)=h_{d-1,1}$.
\endroster
\endproclaim

In [\refMOd] we have defined {\it $MW$-knots} as real algebraic knots
satisfying \CondW. In the same paper we have proved the implication
\CondW~$\Rightarrow$~\CondT.  

As it is seen from the formulation, the implications \CondT$\Rightarrow$\CondPS\ and
\CondT$\Rightarrow$\CondCr\ are purely topological facts. They are specializations
for $(p,q)=(d,d-2)$ of the following results.

\proclaim{ Proposition \propPS }  Let $1\le q\le p$ and $p\equiv q\mod2$. Then $\ps(\Tproj(p,q))=q$.
\endproclaim

This proposition is proven in \S\sectPS.

\proclaim{ Proposition \propMurasugi } {\rm(Murasugi [\refMu])}. Let $1\le q\le p$ and $p\equiv q\mod2$.
Then $\cn(\Tproj(p,q))=p(q-1)/2$.
\endproclaim

\demo{ Proof } By [\refMu; Proposition 7.5],
$\cn(T(p,q))=p(q-1)$. It remains to observe that the double covering of
any projection of $\Tproj(p,q)$ is a projection of $T(p,q)$. \qed
\enddemo

\noindent
{\bf Conjecture \conjMurasugi. } For $1\le q\le p$, any diagram of $T(p,q)$ with $p(q-1)$ crossings
is positive (i.e. all its crossings are positive).
\smallskip

\smallskip
\noindent
{\bf Conjecture \conjRigid. } All $MW$-knots of degree $d$ are rigidly isotopic to $\Tproj(d,d-2)$.

\subhead 1.2. Generalization for arbitrary genus
\endsubhead

Theorem \thMain\ admits a generalizations for algebraic links of any genus.
Given an irreducible algebraic link in $s$ connected components $L=L_1\sqcup\dots\sqcup L_s=\R A$
where $A$ is a {\it dividing} (or {\it Type I\/})
real algebraic curve (this means that $A\setminus \R A$ is not connected; see [\refRo]), it
is convenient to introduce another invariant (see [\refMO]):
$$
      w_\lambda(L) = w(L) + 2\sum_{i<j}\lk(L_i,L_j),                        \eqno(\eqDefWL)
$$
where the linking numbers are considered with respect to a {\it complex orientation} of $\R A$,
that is the boundary orientation induced from one of the two halves into which $\R A$ divides $A$.
In other words, $w_\lambda$ is defined in the same way as $w$, but all the real crossings are counted,
not only those where the both crossing branches belong to the same link component.

For an arbitrary (not necessarily dividing) real algebraic link $L=L_1\sqcup\dots\sqcup L_s=\R A$
we may also define the invariant
$$
      w_{|\lambda|}(L) = |w(L)| + 2\sum_{i<j} |\lk(L_i,L_j)|.
$$
It is clear that the absolute values of all the three invariants $w(L)$, $w_\lambda(L)$, and $w_{|\lambda|}(L)$
are bounded by $N_d-g$ where $g$ is the genus of $A$. Indeed, $N_d-g$ is the maximal possible number
of double points of an irreducible plane projective curve of degree $d$ and genus $g$.

We recall that a non-singular real algebraic curve $A$ is called an $M$-curve if the number of
connected components of $\R A$ attains its maximal possible value $\operatorname{genus}(A)+1$.
It is well-known (and immediately seen from Euler characteristics considerations) that
$M$-curves are dividing.

Theorem~\thAnyGenus\ generalizes Theorem~\thMain\ to the case of arbitrary genus.
In particular, it provides a description of the isotopy type of
links with the maximal (for given degree and genus) value of
$w_\lambda$ and $w_{|\lambda|}$.

For positive integers $a_0,\dots,a_g$ we define $W_g(a_0,\dots,a_g)$ as the (topological) link
$$
    \Big(\sph^3\cap\bigcup_{i=0}^g\big\{(w-c_i z)^{a_i} = \varepsilon z^{a_i+2}\big\}\Big)/(-1)
    \qquad\text{in}\quad \sph^3/(-1)\cong\RP^3
$$
where $c_0,\dots,c_g$ are distinct complex numbers and $0<\varepsilon\ll 1$.
We orient it as the boundary of the complex curve in the 4-ball.
This is a $(g+1)$-component link $K_0\cup\dots\cup K_g$, each component $K_i$ being isotopic
to $\Tproj(a_i+2,a_i)$ and placed in a tubular neighborhood $U_i$
of the line $\rl_i=\{w=c_i z\}/(-1)$ so that $[K_i]=a_i[\rl_i]$ in $H_1(U_i)$
(we orient $\rl_i$ also as the boundary of the complex disk in the 4-ball).
We assume that $\varepsilon$
is so small that the $U_i$'s can be chosen pairwise disjoint.

We call $\rl_0\cup\dots\cup\rl_g$ the {\it Hopf link} in $\RP^3$.
It is indeed a union of fibers of the Hopf fibration
$\eta:\RP^3 = (\C^2\setminus 0)/\R^\times\to(\C^2\setminus0)/\C^\times=\CP^1$.
So, the Hopf link is a union of oriented lines in $\RP^3$ such that any two lines are linked positively
(the linking number is $\frac12$).
We show in Proposition \propHopf\ that it is uniquely determined by this condition.

\proclaim{ Theorem \thAnyGenus }
Let $L$ be an irreducible real algebraic link of degree $d\ge 3$ and genus $g$ which is not contained in a plane.
Then the following conditions are equivalent:
\roster
\item"\CondW"   $w_{|\lambda|}(L)=N_d-g$;
\item"\CondWA"  $\C L$ is dividing and $|w_\lambda(L)| = N_d-g$;
\item"\CondH"   any real plane tangent to $L$ has only real intersections with $\C L$;
\item"\CondPSA" any generic real plane cuts $L$ at $d$ or $d-2$ real points;
\item"\CondPS"  $\ps(L) = d-2$;
\item"\CondCrA" any linear projection of $L$ from a generic point of $\RP^3$
                has at least $N_d-g-1$ hyperbolic double points;
\item"\CondCr"  $\cn(L)=N_d-g-1$;
\item"\CondT"   $L$ or its mirror image is isotopic to
                $W_g(a_0,\dots,a_g)$
                for some positive integers $a_0,\dots,a_g$ such that $a_0+\dots+a_g=d-2$
                (in particular, $g\le d-3$);
\item"\CondTA"  Condition \CondT\ holds in the sense of oriented links for a complex
                orientation of $L$ (note that \CondT\ implies that
                $\C L$ is an $M$-curve, so it has complex orientation);
\endroster
\endproclaim

\smallskip\noindent
{\bf Remark \remWO. } Given any curve $A$ (not necessarily dividing) and any orientation $O$ of $L=\R A$,
the right hand side of (\eqDefWL) is well-defined 
and we denote the resulting
invariant of $(L,O)$ by $w_\lambda(L,O)$.
Then Condition \CondW\ in Theorem \thAnyGenus\ can be replaced (without any change in the proof)
by the condition that $|w_\lambda(L,O)|=N_d-g$ for some orientation $O$ of $L$.
\smallskip

We call irreducible real algebraic links satisfying any of the conditions
\CondW--\CondTA\ of Theorem \thAnyGenus\ {\it $MW_\lambda$-links}. The following corollary shows that,
in a sense, $MW$-links of nonzero genus do not exist.

\proclaim{ Corollary \corAnyGenus } Under the hypothesis of Theorem \thAnyGenus, if $|w(L)|$ attains
the upper bound $N_d-g$, then $g=0$.
\endproclaim

\demo{ Proof } We have $|w|\le w_{|\lambda|}\le N_d-g$. Thus $|w|=N_d-g$ implies \CondW\ and hence \CondTA.
However, if $g>0$, then the linking number of any two components of $W_g(a_0,\dots,a_g)$
is positive, which contradicts the fact that $|w|=w_{|\lambda|}$.
\qed\enddemo

It is interesting to study what is the maximal possible value of $w$ for links of degree $d$ and genus $g$.
The classification obtained in [\refMO]
\footnote{For $d=4$, $g=1$, contrary to an erroneous statement from [9, \S 7.1], there
exists a non-planar link, namely $W_1(1,1)$. Using methods of [\refMO] it is easy to show
that there are exactly 4 rigid isotopy classes: $\pm W_1(1,1)$ ($w_\lambda=\pm2$),
$\bigcirc\bigcirc$ ($w_\lambda=0$), and $\bigcirc$ ($w=0$).}
shows that for $g=1$ and $d=4,5,6$, the maximal value
is attained on $W_1(d-3,1)$ and thus is equal to $(N_d-1)-(d-3)$. It is natural to expect that for any $g\ge0$ and
any $d\ge g+3$, the maximum of $w$ is attained on $W_g(d-g-2,1,\dots,1)$ and thus is equal to
$(N_d-g) - g(d-g-2) - g(g-1)/2$.

\smallskip
The following theorem claims realizability of all isotopy types
allowed by Theorem~\thAnyGenus\ thus finishing the topological classification
of $MW_\lambda$-links of arbitrary degree and genus.

\proclaim{ Theorem \thExist }
For any $d\ge3$ and any positive integers
$a_0,\dots,a_g$ such that $a_0+\dots+a_g=d-2$, there exists an irreducible real algebraic link of degree
$d$ and genus $g$, isotopic to $W_g(a_0,\dots,a_g)$.
\endproclaim

\smallskip
\subhead 1.3. Organization of the paper \endsubhead

We prove Theorem~\thAnyGenus\ in \S\sectProofMainTh\ (Theorem~\thMain\ is, evidently, a partial case of
Theorem~\thAnyGenus). The longest part of the proof concerns the implication
\CondH~$\Rightarrow$~\CondTA\ (\S\sectT). The main idea is to show that if $L$ satisfies \CondH,
then the tangent surface of $L$ (that is the union of all tangent lines) is topologically isotopic
to the boundary of a tubular neighborhood of the Hops link (the isotopy cannot be smooth because
the tangent surface is not: it has a cuspidal edge along $L$; see Figures \figChu\ and \figChuInv).

In \S\sectHopf\ we prove that any configuration of positively linked lines is rigidly isotopic
to the projective Hopf link. This fact is used in the proof of Theorem \thAnyGenus\ but
in our opinion it is interesting by itself.
In \S\sectPermut\ we discuss some topological properties of the projective links $\Tproj(p,q)$.
The subsections \sectDiagr--\sectPermutLast\ are not used in the rest of the paper.
In \S\sectExist\ we prove Theorem \thExist.



\head\sectHopf. Uniqueness of a configuration of positively linked lines in $\RP^3$ up to rigid isotopy
\endhead

In this section a real line is understood in the usual sense, it a subset of $\RP^3$ rather than a
conj-invariant complex line in $\CP^3$.
The following fact is certainly well known for specialists but we did not find any reference.

\proclaim{ Proposition \propHopf } Let $L=\rl_0\cup\dots\cup\rl_n$ and $L'=\rl'_0\cup\dots\cup\rl'_n$
be two unions of $n$ pairwise disjoint oriented lines in $\RP^3$ such that
$\lk(\rl_i,\rl_j)=\lk(\rl'_i,\rl'_j)=\frac12$ for any $i<j$.
Then $L$ is rigidly isotopic to $L'$.
\endproclaim

\noindent
(When speaking of rigid isotopy, we consider $L$ and $L'$ as reducible real algebraic links.)

\demo{ Proof } It is enough to show that $L$ is rigidly isotopic to some standard arrangement of lines.
As a standard arrangement we choose a collection of lines lying on the same hyperboloid.
Note that such an arrangement is a particular case of a Hopf link (see the introduction)
because $\eta^{-1}(C)$ is a hyperboloid when $C$ is a circle on $\CP^1$; recall that
$\eta:\RP^3\to\CP^1$ is the Hopf fibration.

Let us fix an affine chart $(x,y,z)$ so that $\rl_0$ is the line at infinity in the plane \hbox{$x=0$}.
Then (for an appropriate orientation of $\rl_0$) the condition $\lk(\rl_0,\rl_i)>0$ for $i>0$ means that the $x$-coordinate
increases along each of $\rl_1,\dots,\rl_n$, i.e.~all the lines are oriented from the left to the right
(see Figure~\figHopf(a)).
Let us consider the parallel projection
$\pi:(x,y,z)\mapsto(x,y)$.
We assume that the lines are numbered in the descending order according to the
slope of $\pi(\rl_i)$.
Then the positivity of the linking numbers means that each $\rl_i$ passes over $\rl_{1},\dots,\rl_{i-1}$
at the crossings (see Figure~\figHopf(a)).

\midinsert
\centerline{\epsfxsize=90mm\epsfbox{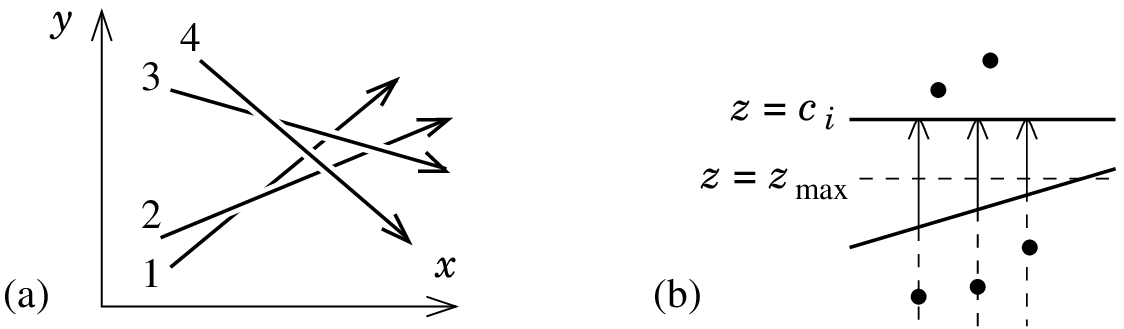}}
\botcaption{ Figure \figHopf } (a) $\pi(L\setminus\rl_0)$; (b) moving $\rl_i$ inside $\pi^{-1}(\rl_i)$
\endcaption
\endinsert

Let $z_{\max}$ be the maximum of the $z$-coordinates of all the preimages of the crossing points,
and let $y=c_i x+ b_i$ be the equation of $\pi(\rl_i)$,
$i=1,\dots,n$. Our assumption about the slopes means that $c_1>c_2>\dots>c_n$.
By shifting vertically all the configuration, we may achieve that $c_n>z_{\max}$.
 
At the first step we move each line $\rl_i$ (one by one starting with $\rl_n$)
inside the plane $\pi^{-1}(\rl_i)$ to place it in the horizontal plane $z=c_i$; see
Figure \figHopf(b) where the lines $\rl_j$ for $j\ne i$ are represented by their
intersection points with the plane $\pi^{-1}(\rl_i)$.
The condition $c_n>z_{\max}$ guarantees that $\rl_i$ does not cross the other lines
during this motion; see Figure \figHopf(b).
Indeed, let $\rl_{ij}=\pi^{-1}(\pi(\rl_i)\cap\pi(\rl_j))$.
Then, during the motion of $\rl_i$, the final 
position of $\rl_i\cap\rl_{ij}$ ($j\ne i$) is over its initial position (because $c_i>z_{\max}$),
and $\rl_j\cap\rl_{ij}$ cannot be between them.

At the second step we parallelly translate each line (staying in the horizontal plane) to make it
cross the axis $x=y=0$. Thus the final position of the $i$-th line is given by the system of two equations
$z=c_i$, $y=c_i x$. This means that all the lines (including $\rl_0$~!) sit on the hyperboloid
$y = zx$.
\qed\enddemo



\head\sectProofMainTh. Classification of $MW_\lambda$-links (proof Theorem \thAnyGenus)
\endhead

\subhead   \sectTors. Non-vanishing of the torsion
\endsubhead

Recall that the sign of the (differential geometric) torsion of a curve
 $t\mapsto r(t)\in\R^3$, $t\in\R$, coincides with the sign of $\det(r',r'',r''')$
and it does not depend on the parametrization if $r'\neq 0$.
The sign of the torsion of a curve in $\RP^3$ does not depend
on a choice of positively oriented affine chart.

\proclaim{ Lemma \lemTors } Let $L=\R A$ be an irreducible real algebraic link of degree $d$ which is not
contained in any plane.
If Condition \CondH\ of Theorem \thAnyGenus\ holds, then the torsion of $L$ nowhere vanishes.
\endproclaim

\demo{ Proof } Let $p\in L$.
We can always choose affine coordinates $(x_1,x_2,x_3)$ centered at $p$ such that $L$ near $p$
admits a parametrization
$$
   t\mapsto(x_1(t),x_2(t),x_3(t)), \quad
      1 = \ord_t x_1 < \ord_t x_2 < \ord_t x_3.
$$
Suppose that the torsion at $p$ is zero, i.e. $a=\ord_t x_3\ge 4$, and let us prove that
then \CondH\ does not hold.

\smallskip
Consider the projection $\pi:(x_1,x_2,x_3)\mapsto(x_1,x_3)$. Then
$\pi(L)$ is parametrized by $t\mapsto\gamma(t)=(t,ct^a+o(t^a))$, $c\ne 0$.
Let $\ell_t$ be the line on the plane $(x_1,x_3)$ tangent to $\pi(L)$ at $\gamma(t)$.
If $0<|t|\ll 1$, them
$\ell_t$ has at least two non-real intersections with $\pi(A)$ near the origin,
thus the plane $\pi^{-1}(\ell_t)$ does not satisfy \CondH.
\qed\enddemo


\subhead \sectHyperbPla. Hyperbolic planar curves
\endsubhead

Let $C$ be a real algebraic curve of degree $d$ in $\RP^2$ (maybe, singular).
Recall that $C$ is called {\it hyperbolic} with respect to a point $q\in\RP^2$
(which may or may not belong to $C$),
if any real line through $q$ intersects $C$ at real points only.
We denote:
$$
   \hyp(C) = \{q\mid\text{$C$ is hyperbolic with respect to $q$}\}.            \eqno(\eqDefHyp)
$$
It is easy to check that $\hyp(C)$ is either empty or a convex closed set. It is possible that $\hyp(C)$
contains only one point. In this case, the point should be singular. For example,
if $C$ is a cuspidal cubic, then $\hyp(C)$ consists of the cusp only.

\proclaim{ Lemma \lemHypC } Let $C$ be a real plane curve, $q\in\hyp(C)$, and
$q_1\in\R C\setminus\{q\}$.
Then each local branch of $C$ at $q_1$ is
smooth, real, and transverse to the line $(qq_1)$. The projection from $q$
defines a covering $\R\tilde C\to\RP^1$ where $\tilde C$ is the normalization of $C$.
\qed
\endproclaim

A real algebraic curve is called {\it dividing} if its normalization is
divided by the real locus into two halves.
The following lemma is due to Rokhlin [\refRo] in the case of smooth curves,
the same proof also holds in the singular case.

\proclaim{ Lemma \lemHypTypeI }
Let $C$ be a real algebraic curve in $\P^2$, and $p\in\hyp(C)$.
Then $C$ is dividing and a complex orientation of $\R C$ is the
pull-back of an orientation of $\RP^1$ under the central projection
$\P^2\setminus\{p\}\to\P^1$ restricted to $\R C$. \qed
\endproclaim


\subhead \sectHyperbSpa. Hyperbolic spatial curves
\endsubhead

Let $L=\R A$ be an irreducible real algebraic link in $\rp^3$ of degree $d$.
The line (in $\P^3$) tangent to $A$ at a point $p$ will be denoted by $T_p$.

Given a point $p\in\RP^3$,
let $\pi_p:\P^3\setminus\{p\}\to\P^2$ be the linear projection from a point $p$ and
let $\hat\pi_p:A\to\P^2$ be the restriction of $\pi_p$ to $A$. 
If $p\in L$, then we extend $\hat\pi_p$ to $p$
by continuity, thus $\pi_p^{-1}(\hat\pi_p(p))=T_p$.

We say that $A$ is {\it hyperbolic} with respect to a real line $\ell$ if, for
any real plane $P$ passing through $\ell$, each intersection point of $A$ and $P\setminus\ell$ is real.
The following property of hyperbolic spatial curves is immediate from the definition:

\proclaim{ Lemma \lemHypKC }
Let $\ell$ be a real line in $\P^3$ and $p\in\R\ell$. Then
$A$ is hyperbolic with respect to $\ell$ if and only if
$\pi_p(\ell)\in\hyp(C_p)$. \qed
\endproclaim

\proclaim{ Lemma \lemHypK }
Suppose that $A$ is hyperbolic with respect to a real line $\ell$ and let
$p\in(\R A)\setminus\ell$. Then $\ell\cap T_p=\varnothing$.
\endproclaim

\demo{ Proof }
Combine Lemma \lemHypC\ and Lemma \lemHypKC.
Namely, let $p_1\in\R\ell$ and $q=\pi_{p_1}(\ell)$. Then $q\in\hyp(C_{p_1})$ by Lemma \lemHypKC.
Hence the line $\pi_{p_1}(T_p)$ cannot pass through $q$ by Lemma \lemHypC\ because
otherwise it would be a line through $q$ non-transverse to the image of the germ $(A,p)$
under the projection $\pi_{p_1}$.
\qed\enddemo

\proclaim{ Lemma \lemHypPosA }
Consider a crossing of a generic projection of $\R A$.
Let $p,q\in\R A$ be the preimages of the crossing point under this projection.
Suppose that $A$ is hyperbolic with respect to the tangent line $T_p A$ and
the torsion of $\R A$ at $p$ is positive.
Then the considered crossing is also positive.
\endproclaim

\midinsert
\centerline{ \epsfxsize=100mm\epsfbox{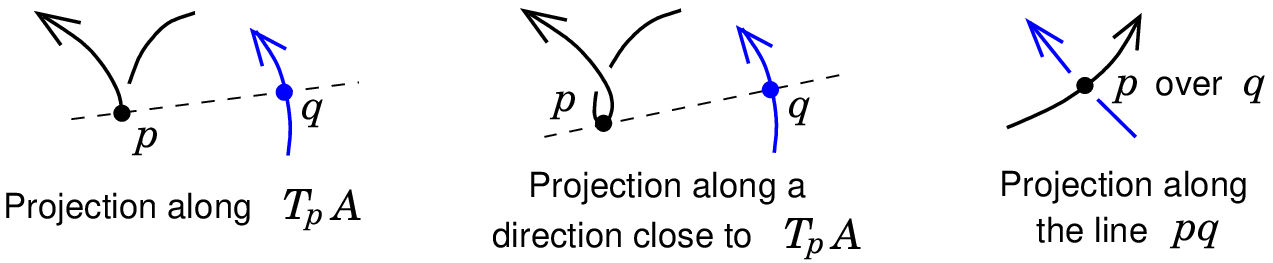}}
\botcaption{ Figure \figHypPos } Positive crossing and positive torsion at $p$
\endcaption
\endinsert

\demo{ Proof }
Briefly speaking, the statement follows from Lemma \lemHypTypeI\ and Figure \figHypPos.
Indeed, since the torsion at $p$ is positive,
the curve near $p$ looks as in Figure \figHypPos~(left)
because a small perturbation of this projection should produce a positive crossing near $p$
as shown in Figure \figHypPos~(middle).
Since the curve is hyperbolic with respect to $T_p A$,
by Lemma~\lemHypTypeI\ the mutual orientations at $p$ and $q$ should be as in Figure \figHypPos~(left),
i.e., when looking from $p$ to $q$, the rotation from $\vec v_p$ to $\vec v_q$ is counterclockwise
where $\vec v_p$ and $\vec v_q$ are positive tangent vectors at $p$ and $q$ respectively.
This means exactly that the crossing is positive; see Figure \figHypPos~(right).
\qed\enddemo

\proclaim{ Lemma \lemHypPos }
Let $\R A$ be an irreducible real algebraic link in $\RP^3$ which is
hyperbolic with respect to any real tangent line (and thus dividing by
by Lemma \lemHypTypeI\ and Lemma \lemHypKC).
Suppose that the torsion of $\R A$ is nonzero at any point.
Then the torsion of $\R A$ at all its points has the same sign.
Furthermore all non-isolated crossings of any
projection have the same sign as the torsion of $\R A$ once we endow
$\R A$ with a complex orientation.

\endproclaim

\demo{ Proof } For any two generic points $p,q\in\R A$, the torsion at these points has the same sign
due to Lemma \lemHypPosA\ applied to the projection along the line $pq$.
Hence the torsion is positive everywhere. Then the positivity at any crossing of any projection
follows again from Lemma \lemHypPosA.
\qed\enddemo


\subhead \sectT. Tangent surface
\endsubhead
In this subsection we prove the implication \CondH~$\Rightarrow$~\CondTA\ in  Theorem \thAnyGenus.
Let the notation be as introduced in the beginning of \S\sectHyperbSpa\
and assume that Condition \CondH\ of Theorem \thAnyGenus\
holds, i.e., any real tangent plane has only real intersections with $A$, in particular,
$A$ is hyperbolic with respect to $T_p$ for any $p\in L=\R A$.
By Lemma \lemHypTypeI\ combined with Lemma \lemHypKC, the curve $A$ is dividing, so we
fix a complex orientation on $L$ and we orient the real locus of each tangent line $\R T_p$, $p\in L$,
accordingly.

The torsion of $L$ does not vanish by Lemma \lemTors, thus
lemma \lemHypPos\ implies that the torsion if $L$ has the same sign everywhere.
We assume that it is positive (otherwise we replace $L$ by its mirror image).

\proclaim{ Lemma \lemHypII }
Let $p\in L$. Then $A\cap T_p = \{p\}$.
\endproclaim

\demo{ Proof }
Suppose that there exists a point $q\in (A\cap T_p)\setminus\{p\}$.
Then $q$ is real by \CondH\ applied to a tangent plane at $p$,
and $A$ is hyperbolic with respect to $T_q$ by \CondH\ applied to a tangent plane at $q$.
If $T_p\ne T_q$, then $p\not\in T_q$, and we obtain a contradiction with Lemma \lemHypK\ for $\ell=T_q$.
Hence $T_p=T_q$.

Let $\pi$ be a generic projection. Then $\pi(A)$ is tangent to $\pi(T_p)$ at $\pi(p)$ and $\pi(q)$.
Let $\ell$ be a line tangent to $\pi(T_p)$ at a point close to $\pi(p)$.
By Lemma \lemTors\ we may assume that the curvature of $\pi(A)$ at $\pi(p)$ and $\pi(q)$ is non-zero,
hence $\ell$ can be chosen so that it has two non-real intersections with $\pi(A)$ near $\pi(q)$.
Then the plane $\pi^{-1}(\ell)$ contradicts \CondH. 
\qed\enddemo

\proclaim{ Lemma \lemHypIII }
Let $p\in L$, $q\in A\setminus L$, and let $\ell$ be the line $(q\bar q)$.
Then $T_p\cap\ell = \varnothing$.
\endproclaim

\demo{ Proof } If $T_p\cap\ell \ne \varnothing$, then \CondH\
does not hold for the plane containing these two lines.
\qed\enddemo

\proclaim{ Lemma \lemDisjoint } Let $p$ and $q$ be two distinct points on $L$.
Then $T_{p}\cap T_{q}=\varnothing$
and $\lk(\R T_p,\R T_q)=1/2$.
\endproclaim

\demo{ Proof } Let $\ell=T_q$. Then $A$ is hyperbolic with respect to $\ell$ by Condition \CondH, and
$p\not\in\ell$ by Lemma \lemHypII.
Thus $T_p\cap\ell=\varnothing$ by Lemma \lemHypK.
The positivity of the linking number follows from Lemma \lemHypPos.
\qed\enddemo

\proclaim{ Lemma \lemTwoLines } For any $p\in L$
there exist two real lines $\ell$ and $\ell'$ arbitrarily close to $T_p$
such that $A$ is hyperbolic with respect to each of them, $\ell\cap A=\varnothing$, and
$\ell'$ crosses $A$ without tangency at a pair of complex conjugated points.
\endproclaim

\demo{ Proof }
Let $p_0\in\R T_p\setminus\{p\}$.
The curve $C_{p_0}$ has an ordinary cusp at $q_0=\pi_{p_0}(p)$ because $L$ has non-zero torsion at $p$.
Let $p_1$ and $p_2$ be generic points close to $p_0$ which are chosen on different sides of
the osculating plane of $L$ at $p$. Then $C_{p_1}$ and $C_{p_2}$ are
obtained from $C_{p_0}$ by a perturbation of the cusp as shown in Figure \figCusp\ where
$q_2$ is a solitary node of $C_{p_2}$ (a point where two conjugate non-real local branches cross).

\midinsert
\epsfxsize=60mm
\centerline{\epsfbox{ 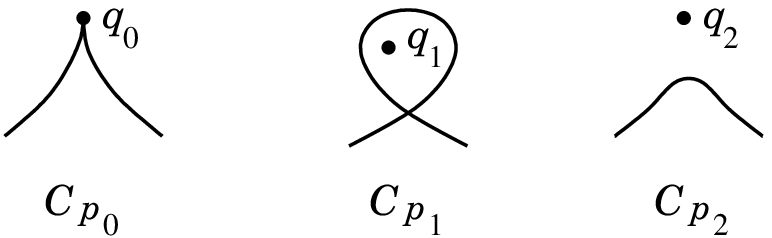 }}
\botcaption{ Figure \figCusp } Two perturbations of a cusp in the proof of Lemma \lemTwoLines
\endcaption
\endinsert

Let $\ell = \pi_{p_1}^{-1}(q_1)$ and $\ell' = \pi_{p_2}^{-1}(q_2)$
where the points $q_1$ and $q_2$ are chosen as in Figure \figCusp.
The fact that $A$ is hyperbolic with respect to $T_p$ implies $q_0\in\hyp(C_{p_0})$ by Lemma \lemHypKC.
We see in Figure \figCusp\ that then $q_j\in\hyp(C_{p_j})$ for $j=1,2$.
\qed\enddemo

Let $K_0,\dots,K_n$ be the connected components of $L$.
For each  $i=0,\dots,n$, let $\ell_i$ and $\ell'_i$  be the lines given by Lemma \lemTwoLines\
close to $T_p$ for some $p\in K_i$.
To simplify the notation, we denote $\R\ell_i$ and $\R\ell'_i$ by $\rl_i$ and $\rl'_i$.
We endow $\rl_i$ and $\rl'_i$ with the orientation inherited from $\R T_p$ and we set
$$
        a_i=2\lk(\rl'_i,K_i).
$$

Let
$$
       T=TK_0\cup\dots\cup TK_n \quad\text{ where }\quad
       TK_i = \bigcup_{p\in K_i} \R T_p\,.
$$
By Lemma \lemDisjoint, $T$ is a disjoint union of a continuous family
of real projective lines (topologically, circles) parametrized by $L$.
We are going to show that $T$  is topologically isotopic in $\RP^3$
to a union of hyperboloids.
Note that $T$ is not smooth. It has a cuspidal edge along $L$
(see Figures \figChu--\figTwoChu\ and Remark \remCycloid).

\proclaim{ Lemma \lemLK }
For any $i$ we have $a_i>0$ and
$$
\xalignat3
     &2\lk(\rl'_i,L  ) = d-2,   & &2\lk(\rl_i,L)=d,\\
     &2\lk(\rl'_i,K_i) = a_i,   & &2\lk(\rl_i,K_i)=a_i+2,\\
     &2\lk(\rl'_j,K_i) = a_i,   & &2\lk(\rl_j,K_i)=a_i       &&\text{for $j \ne i$.}
\endxalignat
$$
\endproclaim

\demo{ Proof }
The positivity of $a_i$ follows from the hyperbolicity of $L$ with respect to $\ell'_i$.
The first two rows of identities in Lemma \lemLK\ follow from hyperbolicity
and the definition of $a_i$.
We also have $\lk(\rl_j,K_i)=\lk(\rl'_j,K_i)$ for $j\ne i$ because both $\rl_j$, $\rl'_j$ are
close to $T_p$ for some $p\in K_j$, and $T_p$ is disjoint from $K_i$.

So, it remains to show that $|\lk(\rl'_i,K_i)|=|\lk(\rl'_j,K_i)|$.
Indeed, each of $\ell_i'$, $\ell'_j$ passes through a conjugate pair of points of $A$.
Since $A\setminus L$ has two connected components, there exists a continuous path
$p:[0,1]\to A\setminus L$ such that $\rl'_i=\rl(0)$ and $\rl'_j=\rl(1)$ where
$\rl(t)$ is the real locus of the line $\big(p(t)\,\overline{p(t)}\big)$. Each line $\rl(t)$
is disjoint from $L$ by Lemma \lemHypIII, hence $|\lk(\rl(t),K_i)|$ is constant.
\qed\enddemo

By Lemma \lemLK\ we have
$$
    a_0+\dots+a_n
    = \sum_{i=0}^n\lk(\rl'_1,K_i)=\lk(\rl'_1,L)=d-2.         \eqno(\eqSumA)
$$

Let $i\in\{0,\dots,n\}$.
The line $\rl'_i$ is disjoint from $T$ by Lemma \lemHypK.
Let $P$ be a real plane through $\rl'_i$. 
Again by Lemma \lemHypK, $P$ crosses each line $T_p$, $p\in K_i$, at a single point.
Let us denote this point by $\xi_{P}(p)$. Then $\xi_{P}:K_i\to\R P$ is a continuous mapping.
It is injective by Lemma \lemDisjoint\ and its image (which is $TK_i\cap\R P$) is disjoint from $\rl'_i$.
Hence $TK_i\cap\R P$ is a Jordan curve in the affine real plane $\R P\setminus\rl'_i$.
Let $D_{P}$ be the disk bounded by this Jordan curve, and let $U_i=\bigcup_P D_{P}$ where
$P$ runs through all the real planes through $\rl'_i$.
Then $U_i$ is fibered by disks over a circle which parametrizes the pencil of planes through $\rl'_i$.
Since $\RP^3$ is orientable, this fibration is trivial, thus $U_i$ is a solid torus and $\partial U_i=T(K_i)$.
Each $P$ transversally crosses $K_i$ at $a_i$ real points, thus $K_i$ sits in $TK_i$ and
it realizes the homology class $a_i\alpha_i$ where $\alpha_i$ is a generator of $H_1(U_i)$.

The same arguments applied to the line $\rl_i$ show that $TK_i$ bounds a solid torus $V_i$
such that  $K_i$ realizes the homology class $(a_i+2)\beta_i$ where $\beta_i$ is a generator of $H_1(U'_i)$.
Easy to see that $\alpha_i=[\rl_i]\in H_1(U_i)$ and $\beta_i=[\rl'_i]\in H_1(V_i)$.
Using Proposition \propHabTpq\ we deduce that $K_i$ realizes the knot $\Tproj(a_i,a_i+2)$
sitting in the torus $TK_i$.

By comparing the linking numbers given by Lemma \lemLK\ (or using the fact that
all the $\rl'_i$ are homotopic to each other in the complement of $T$; see the proof of Lemma \lemLK),
we see that the solid tori $U_0,\dots,U_n$ are pairwise disjoint, so their union can be considered as
a tubular neighborhood of the line arrangement $\Cal L = \rl_0\cup\dots\cup\rl_n$
(see Figure \figTwoChu).
All the linking numbers $\lk(\rl_i,\rl_j)$, $i\ne j$, are positive by Lemma \lemHypPos,
hence $\Cal L$ is the Hopf link by Proposition \propHopf.

\midinsert
\centerline{
   \epsfxsize=40mm\epsfbox{ 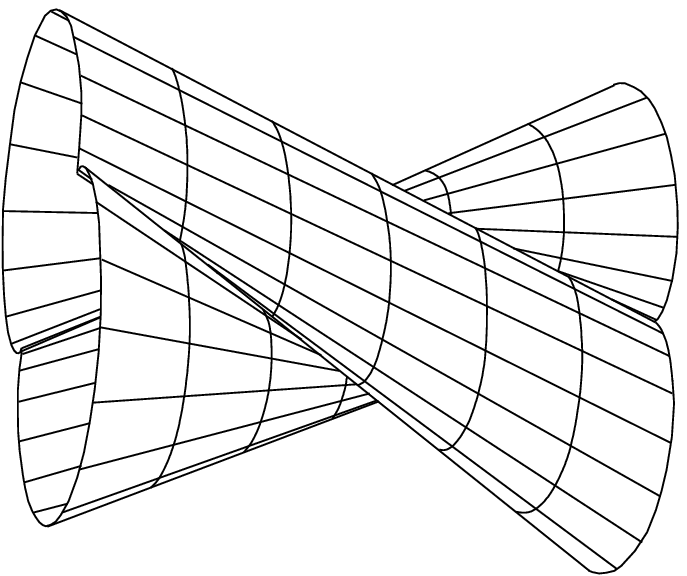 }\hskip 10mm
   \epsfxsize=45mm\epsfbox{ 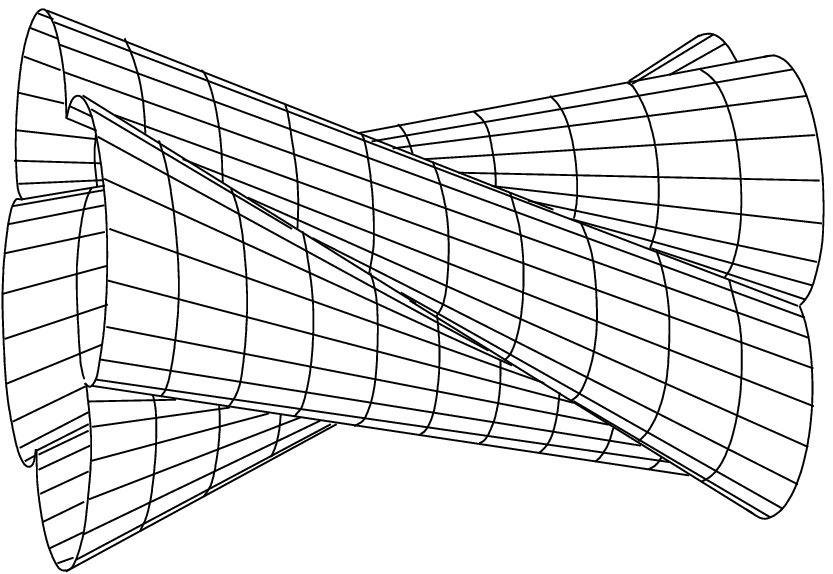 }
}
\botcaption{ Figure \figChu } Tangent surface of $\Tproj(4,2)$ and $\Tproj(6,4)$ (view from $V$) 
\endcaption
\endinsert

\midinsert
\centerline{
   \epsfxsize=50mm\epsfbox{ 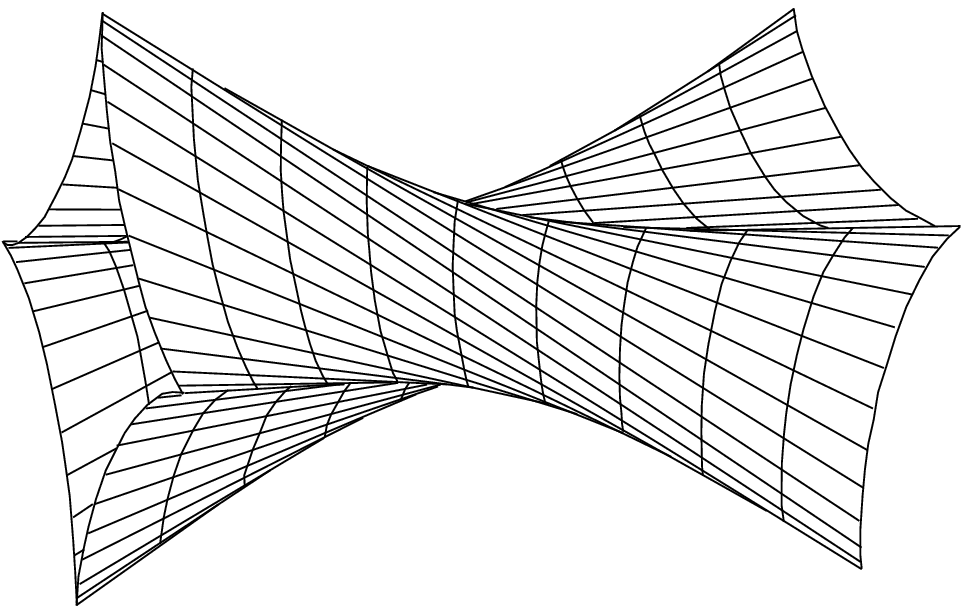 }\hskip 3mm
   \epsfxsize=65mm\epsfbox{ 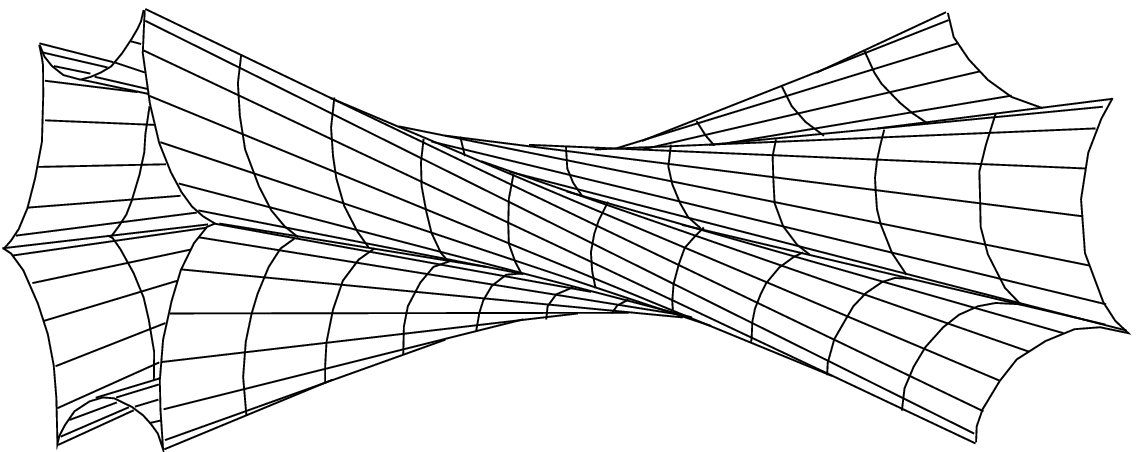 }
}
\botcaption{ Figure \figChuInv } Tangent surface of $\Tproj(4,2)$ and $\Tproj(6,4)$ (view from $U$) 
\endcaption
\endinsert

\midinsert
\centerline{ \epsfxsize=50mm\epsfbox{ 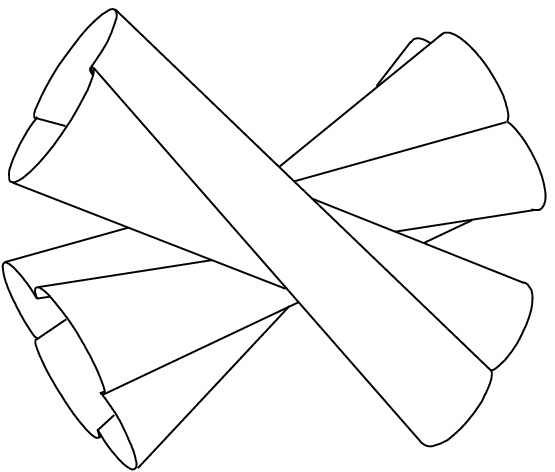 } }
\botcaption{ Figure \figTwoChu } Mutual position of $TK_0$ and $TK_1$ for $W_1(4,2)$
\endcaption
\endinsert

It remains to show that $n=g=\text{genus}(A)$. Indeed, we have
$$
\xalignat2
   \cn(L)&\ge\sum_{i=0}^n\cn(K_i) + 2\sum_{i<j}\lk(K_i,K_j) \\
           &=\tfrac12\sum_{i=0}^n{(a_i+2)(a_i-1)} + \sum_{i<j}a_i a_j &&\text{by Proposition \propMurasugi}\\
           &=\tfrac12\Big(\sum_{i=0}^n a_i\Big)^2 + \tfrac12\sum_{i=0}^n a_i - (n+1)\\
           &= N_d - n - 1  && \text{by (\eqSumA)}
\endxalignat
$$
(note that ``$\ge$" can be replaced here by ``$=$" because the minimal number of crossings
for each $K_i$ can be realized inside $U_i$).
On the other hand, the number of double points of any projection is bounded by $N_d-g$ by genus formula,
and there exists a projection with at least one solitary double point (it is enough to project from
a real point on a line $(p\bar p)$ for $p\in A\setminus\R A$), hence $\cn(L)\le N_d-g-1$.
Thus
$n\ge g$. Since $n\le g$ by Harnack inequality, we conclude that $n=g$, hence
$L\cong W_g(a_0,\dots,a_g)$.
The implication \CondH\ $\Rightarrow$ \CondTA\ is proven.

The above computation of $\cn(L)$ also shows that
$$
    \cn(W_g(a_0,\dots,a_g)) = N_d-g-1                            \eqno(\eqCRW)
$$
for any $g\ge 0$ and any positive $a_i$'s satisfying $a_0+\dots+a_g=d-2$.


\subhead \sectEOP. End of the proof of Theorem \thAnyGenus
\endsubhead
In Figure \figSchPrf\ we show the implications which are either
evident or already proven. Here we prove the others.

\midinsert
\centerline{\epsfxsize=65mm\epsfbox{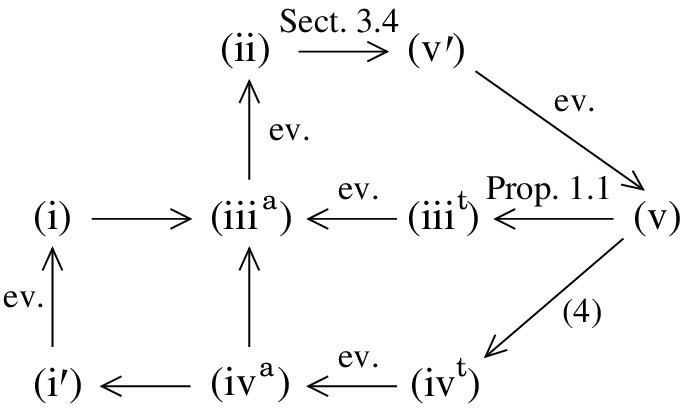}}
\botcaption{ Figure \figSchPrf } Scheme of implications (``ev.'' means ``evident'')
\endcaption
\endinsert

\medskip

\CondW\ or \CondCRA\ $\Rightarrow$ \CondPSA.
Suppose that \CondPSA\ is not satisfied. Let $P$ be a generic real plane which cuts $L=\R A$
at less than $d-2$ real points.
Then there are at least two conjugate pairs of non-real points $p, \bar p$
and $q,\bar q$ in $P\cap A$.
Since $P$ is generic, we may assume that it does not contain any trisecant of $A$.
Then the lines $p q$ and $\bar p\bar q$ are distinct.
Since they are contained in $P$, they cross each other.
Then the projection from their intersection point has at least two
non-real double points. This contradicts each of the conditions \CondW\ and \CondCRA.
\medskip

\medskip
\CondH\ $\&$ \CondCRA\ $\Rightarrow$ \CondWA.
Recall that a double point (node) of a real curve in the plane is called {\it hyperbolic}
if it is an intersection point of two real branches.
A projection of $A$ from a smooth point $p$ on the tangent surface $T$ has a cusp.
If we push $p$ from $T$ then the cusp in the resulting projection gets replaced
with a double point that is either hyperbolic or solitary depending on the choice of
the side for moving $p$ (see Figure \figCusp). By \CondCRA\ the number of hyperbolic nodes
is at least $N_g-g-1$.
Thus this number is $N_g-g$ for
one of the choices.
By Lemma \lemHypPos\ (using the condition \CondH\ here as well
as for deducing that $A$ is dividing) all these nodes must correspond
to positive crossing points of the diagram.



\subhead\sectPropU. Concluding remarks for Theorem \thAnyGenus
\endsubhead

We see that $T$ cuts $\RP^3$ into a union of solid tori $U=U_0\cup\dots\cup U_g$ and
its complement $V=V_0\cap\dots\cap V_g$.

\proclaim{ Proposition \propU } {\rm(Compare with Condition \CondCRA).}
Let $p$ be a generic point of $\RP^3$. Then $C_p$ has only real double points.
If $p\in U$, then all the double points are hyperbolic 
and the interior of $\hyp(C_p)$ is non-empty.
If $p\in V$, then one double point $q$ is solitary (i.e. has a conjugate pair of local branches),
all the other double points are hyperbolic,
and $\hyp(C_p)=\{q\}$.
\endproclaim

\demo{ Proof } Let us consider a generic path $p(t)$ which connects the given point to a point on $T$.
It defines a continuous deformation of the knot diagram which is a sequence of Reidemeister moves
(R1) -- (R3). However, (R2) is impossible because it involves a negative crossing, and (R1) may
occur only when $p(t)$ passes through $T$. Thus the number and the nature of double points
does not change during the deformation. The projection from a point of $T$ is cuspidal and it is
hyperbolic with respect to the cusp, so all the double points are hyperbolic by Lemma \lemHypC.

Non-emptiness of the interior of $\hyp(C_p)$
in case $p\in U$, follows from the fact that $\hyp(C_p)$
can disappear only by a move (R3). This is however impossible because all crossings are positive
and the boundary orientation on $\partial(\hyp(C_p))$ agrees with an orientation of $C_p$
due to Lemma \lemHypC\ (see Figure \figR).

The hyperbolicity of $C_{p(t)}$ with respect to the solitary node cannot fail
during the deformation due to Lemma \lemHypIII.
\qed\enddemo

\midinsert
\centerline{\lower-10mm\hbox{$\hyp(C_p)$}\epsfxsize=20mm\epsfbox{ 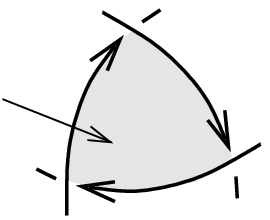 }}
\botcaption{ Figure \figR } Impossibility of a move (R3) which eliminates $\hyp(C_p)$
\endcaption
\endinsert

The sane arguments prove the following fact.

\proclaim{ Proposition \propFibr }
The family of real lines $(p\bar p)$, $p\in A\setminus L$, defines
a fibration of $V$ over the quotient space $(A\setminus L)/\conj$. \qed
\endproclaim

\midinsert
\centerline{\epsfxsize=50mm\epsfbox{ hypocyc.eps } \hskip 10mm
            \epsfxsize=50mm\epsfbox{ epicyc.eps }}
\botcaption{ Figure \figCycloid } $n$-hypocycloid (left) and $n$-epicycloid (right) for $n=5$
\endcaption
\endinsert

\smallskip\noindent
{\bf Remark \remCycloid.} It is a nice and easy exercise 
to check that, in the case of the most symmetric (rational) $MW$-knot,
the Jordan curves $T\cap P$ (discussed after Lemma \lemLK)
are $d$-hypocycloids (for $P\supset \rl'_1$)
or $(d-2)$-epicycloids (for $P\supset\rl_1$).
These are curves defined as trajectories of a point of a smaller circle which
rolls without slipping around a bigger one (see Figure \figCycloid).
By ``the most symmetric $MW$-knot'' we mean $K=\{z^d=w^{d-2}\}$ in $\RP^3=\sph^3/(-1)$,
and $\rl'_1=\{z=0\}$, $\rl_1=\{w=0\}$.
In Figures \figChu\ and \figChuInv\ we show the true shape of the tangent surface $T$ for this knot
for $d=4$ and $d=6$.
The grid in these figures is composed of the tangent lines $T_p$, $p\in K$, and
the plane sections which are hypo- or epicycloids.



\head\sectPermut. Some topological properties and characterizations of $\Tproj(d,d-2)$
\endhead

\subhead\sectHabTpq. Equivalent definitions of $\Tproj(p,q)=h_{a,b}$
\endsubhead

Let $H$ be a hyperboloid in $\RP^3$. It cuts $\RP^3$ into two solid tori $U$ and $V$.
Let $u$ and $v$ be core circles of $U$ and $V$ respectively. We orient them so that
$\lk(u,v)=1/2$. The surface $H$ is doubly ruled.
The two families of projective lines on $H$ give a diffeomorphism
$H\approx \RP^1\times \RP^1$ as well as a basis $(\alpha,\beta)$ of $H_1(H;\Z)$ (defined up to transposition
of basis vectors and changing their sign). We choose it so that any two lines representing $\alpha$
are linked positively
while $\lk(\alpha,u)=\frac12$, $\lk(\beta,u)=\frac12$, $\lk(\alpha,v)=\frac12$, $\lk(\beta,v)=-\frac12$.
Note that the definition of the link $\Tproj(p,q)$ given in the introduction extends
to the case when $p$ and $q$ are of arbitrary signs (but, however, $(p,q)\ne(0,0)$)
by replacing $z^p$ with $\bar z^{|p|}$ (resp. $w^q$ with $\bar w^{|q|}$) if $p<0$ (resp. if $q<0$).
We have isotopies: $\Tproj(p,q)\sim\Tproj(-p,-q)\sim\Tproj(q,p)$. 
The following fact is very easy and we omit its proof (see some more details in [\refMO; \S3]).

\proclaim{ Proposition \propHabTpq } Let $L$ be a smooth oriented one-dimensional submanifold of $H$.
Let $p=a+b$, $q=a-b$ where $a,b\in\Z$, $(a,b)\ne (0,0)$.
We assume that any non-empty one-dimensional submanifold of $L$ is non-trivial in $H_1(H;\Z)$.
If any of the following conditions holds,
then $L$ is isotopic to $\Tproj(p,q)=h_{a,b}$:
\roster
\item"(i)"   $[L]=a\alpha+b\beta$ in $H_1(H;\Z)$;
\item"(ii)"  $\lk(L,u)=p/2$ and $\lk(L,v)=q/2$;
\item"(iii)" $[L]=q[u]$ in $H_1(U;\Z)$ and $[L]=p[v]$ in $H_1(V;\Z)$. \qed
\endroster
\endproclaim

Any of the properties (i)--(iii) can be chosen for a definition of $\Tproj(p,q)$
considered up to isotopy.

\subhead\sectPS. Plane section number of a projective torus link in $\RP^3$
\endsubhead
The following is an equivalent statement of Proposition \propPS.

\proclaim{ Proposition \propSect } The plane section number of $\Tproj(p,q)$
is equal to $\min(|p|,|q|)$.
\endproclaim

\demo{ Proof } Let the notation be as in Section \sectTpqHab\ and let
$L$ be as in Proposition \propTpqHab. Let $P$ be a generic surface in $\RP^3$
isotopic to a plane $\rp^2\subset\rp^3$.
Then either $P\cap U$ or $P\cap V$ is orientable
(since each non-orientable component of $P\setminus H$
contains a M\" obius band whose complement is a disk), let it be $P\cap U$.
Let $C=P\cap H$ oriented as the boundary of $P\cap U$.
Let $j_*:H_1(H)\to H_1(U)$ be the homomorphism induced by the inclusion.
Then $[C]\in\ker j_*$ and $\ker j_*$ is generated by $\alpha+\eps\beta$, $\eps=\pm1$,
thus $[C]=k(\alpha+\eps\beta)$, $k\in\Z$.

Note that $k$ must be odd since the $\Z_2$-homology class
of a plane section of $H\subset\RP^3$ is defined invariantly
of $P$ and thus $[C]\equiv \alpha + \beta\mod 2$.
In particular, $k\neq 0$. By Proposition \propHabTpq\ we have
$[L]=\frac{p+q}2\alpha+\frac{p-q}2\beta$.
Therefore $|[L].[C]|=|k(\frac{p+q}2\pm\frac{p-q}2)|\ge\min(|p|,|q|)$.
\qed\enddemo

\subhead\sectDiagr. Diagrams of projective links and projective braid closures
\endsubhead
Similarly to links in $\sph^3$, a link in $\RP^3$ can be represented by a diagram in a (round) disk $\Bbb D$.
It is a union of immersed circles and/or arcs with endpoints on $\partial\Bbb D$
transverse to each other. The crossings are presented by over- and underpasses
in the same way as for the usual links. The endpoints are placed on $\partial\Bbb D$
symmetrically with respect to the center. 
The link in $\RP^3$ represented by a diagram is defined as follows.
First we consider a tangle in the 3-ball represented by the diagram assuming that
the endpoints of the tangle are placed on the equator. Then we identify the opposite points
of the 3-ball's boundary, see [\refD] for more details.

We define the braid closure in $\sph^3$ in the usual way and we define the braid closure in $\RP^3$
(the {\it projective braid closure\/}) as the link in $\RP^3$ represented by a diagram of
the braid placed in a disk with the endpoints of the strands on the boundary.

Let $p$ and $q$ be positive integers of the same parity.
Then $T(p,q)$ is the closure in $\sph^3$
of the $p$-braid $(\alpha\beta)^q$ where
$\alpha=\sigma_1\sigma_3\dots$ and $\beta=\sigma_2\sigma_4\dots$.
Similarly, $\Tproj(p,q)$ is the closure in $\RP^3$ of the braid represented by the first half of the
word $(\alpha\beta)^q$, i.~e., the braid
$$
   t_{p,q} =
   \cases (\alpha\beta)^{q/2}            &\text{if $p$ and $q$ are even,}\\
          (\alpha\beta)^{(q-1)/2}\alpha  &\text{if $p$ and $q$ are odd,}
   \endcases                                                                   \eqno(\eqT)
$$
see Figure \figT. We obtain this braid by projecting $\Tproj(p,q)$ (defined by the formulas in the
introduction) from a point on a suitable coordinate axis.

\midinsert
\centerline{\hskip 0pt
  \epsfxsize=22mm\hbox to 26mm{\hfill\epsfbox{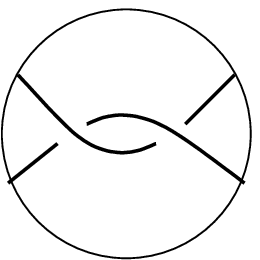}\hfill}\hskip5mm
  \epsfxsize=22mm\hbox to 26mm{\hfill\epsfbox{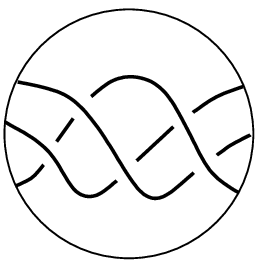}\hfill}\hskip5mm
  \epsfxsize=22mm\hbox to 26mm{\hfill\epsfbox{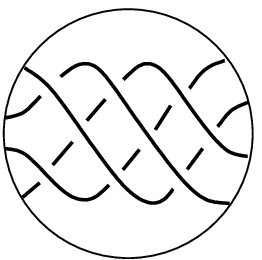}\hfill}\hskip5mm
  \epsfxsize=22mm\hbox to 26mm{\hfill\epsfbox{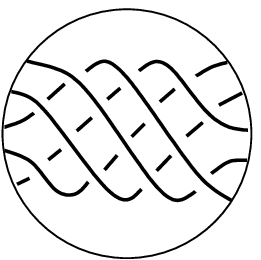}\hfill}}
\centerline{\hskip 0pt
  \hbox to 26mm{\hfill $t_{2,4}=(\alpha\beta)^2$\hfill}\hskip5mm
  \hbox to 26mm{\hfill $t_{3,5}=(\alpha\beta)^2\alpha$\hfill}\hskip5mm
  \hbox to 26mm{\hfill $t_{4,6}=(\alpha\beta)^3$\hfill}\hskip5mm
  \hbox to 26mm{\hfill $t_{5,7}=(\alpha\beta)^3\alpha$\hfill}}
\centerline{\hskip 0pt
  \hbox to 26mm{\hfill $(\alpha,\beta)=(\sigma_1,1)$\hfill}\hskip5mm
  \hbox to 26mm{\hfill $(\sigma_1,\,\sigma_2)$\hfill}\hskip5mm
  \hbox to 26mm{\hfill $(\sigma_1\sigma_3,\,\sigma_2)$\hfill}\hskip5mm
  \hbox to 26mm{\hfill $(\sigma_1\sigma_3,\,\sigma_2\sigma_4)$\hfill}}
\botcaption{ Figure \figT } Knots $\Tproj(d,d-2)$ as projective braid closures
\endcaption
\endinsert

\subhead\sectPermut.4.  Permutation braids
\endsubhead
Let $B_n$ be the braid group on $n$ strands. Let $\mu=\mu_n:B_n\to S_n$ be its
standard homomorphism to the symmetric group (which takes each $\sigma_i$ to the transposition $(i,i+1)$)
and let $e=e_n:B_n\to\Z$ be the {\it exponent sum} homomorphism (which takes each $\sigma_i$ to 1).

{\it Permutation braids} are introduced by El Rifai and Morton in [\refEM]. They can be defined
as braids represented by diagrams whose any two strands cross at most
once and all crossings are positive. One can choose such a diagram so that
each strand is represented by a  straight line segment.
Note that $\mu_n$ establishes a bijection of the set of permutation $n$-braids
onto $S_n$;
see more details in [\refEM; Section 2].

The permutation $n$-braid corresponding to the permutation $i\mapsto n+1-i$, $i=1,\dots,n$, is a positive half-twist
of all the strands
$\Delta=\Delta_n=\prod_{i=1}^{n-1}\prod_{j=1}^{n-i}\sigma_j$ called also the {\it Garside element}
of the braid group $B_n$. It is a unique permutation braid whose exponent sum attains the
maximal possible value $n(n-1)/2$.

\subhead\sectPermut.5. Links in $\RP^3$ and their lifts in $\sph^3$
\endsubhead
Let $L$ be a link in $\RP^3$ and let $\tilde L$ be its inverse image under the covering
$\sph^3\to\RP^3$.

\medskip\centerline
 {\sl Does $\tilde L$ determine $L$?}

\medskip\noindent
Surprisingly, this question seems to be open.
The same question can be posed about links in lens spaces $L(p,q)$
and their lifts to $\sph^3$ (note that $L(2,1)$ is $\RP^3$) but,
as shown in [\refMa], the answer is negative for some knots in $L(4,1)$
and in $L\big(p,{p\pm1\over2}\big)$ for odd $p\ge 5$.
So, the question is not as simple as one could expect.

However, due to a result by Gonz\'alez-Meneses [\refGM],
the answer in a sense is positive if we speak of braids instead of links.
Namely, let $\tau:B_n\to B_n$ be the inner automorphism
$X\mapsto\Delta X\Delta^{-1}$.
We have $\tau(\sigma_i)=\sigma_{n-i}$, $1\le i<n$, thus
$\tau(X)$ is $X$ rotated by $180^\circ$.
We write $X\sim Y$ if the braids $X$ and $Y$ are conjugate.

\proclaim{ Proposition \propLift }
(a). The conjugacy class of $X\Delta$ determines the projective braid closure of
a braid $X$.

\smallskip
(b). The lift of the projective braid closure of $X$
is the braid closure of $X\,\tau(X)$.

\smallskip
(c). The conjugacy class of $X\tau(X)$ determines the closure of $X$ in $\RP^3$.
\endproclaim

\demo{Proof}
(a,b). Evident.

\smallskip
(c). Since $\Delta^2$ belongs to the center of $B_n$, we have
$$
    X\,\tau(X)\sim Y\tau(Y) \;\Leftrightarrow\; (X\Delta)^2\sim(Y\Delta)^2,
$$
thus the result follows from Gonz\'alez-Meneses' theorem [\refGM] which states that
the $m$-th root of a braid (for any $m$) is unique up to conjugacy.
\qed\enddemo

\subhead\sectPermut.6. Characterizations of $\Tproj(d,d-2)$ via braids
\endsubhead
A braid is called {\it positive} if it can be represented by a diagram with positive crossings only
(i.e. if it is a product of the standard generators in positive powers).
We denote the monoid of positive $n$-braids by $B_n^+$.
Recall that $N_d = (d-1)(d-2)/2$.

\proclaim{ Theorem \thPermut }
(a). Let $X$ be a permutation $d$-braid such that $e(X)=N_d$ and the projective closure of
$X$ is a knot {\rm(the last condition is equivalent to the fact that $\mu_d(X\Delta_d)$ is a $d$-cycle).}
Then the projective closure of $X$ is $\Tproj(d,d-2)$.

\smallskip
(b). Let $X$ be a $(d-2)$-braid such that $e(X)=N_d-1$, $X=\Delta A$ for some $A\in B^+_{d-2}$,
and the projective closure of $X$ is a knot.
Then the projective closure of $X$ is $\Tproj(d,d-2)$.
\endproclaim

\noindent{\bf Remark \remPermut.} Using the algorithm in [\refEM],  for any given braid $X$
it is very easy to check if
the condition $X=\Delta A$, $A\in B^+_{d-2}$, in Theorem \thPermut(b) is fulfilled or not.

\smallskip\noindent{\bf Remark \remNonPermut.}
A positive $d$-braid whose projective closure is $\Tproj(d,d-2)$ need not be
a permutation braid. Example: $\sigma_{d-2}^{-1}\Delta_{d-1}\sigma_2$ viewed as a $d$-braid, $d\ge4$.

\proclaim{ Lemma \lemCoxeter } Let $x_1,\dots,x_n$ be elements of a group $G$
such that $x_i x_j = x_j x_i$ for $j-i \ge 2$.
Let $(y_1,\dots,y_n)$ be a permutation of $(x_1,\dots,x_n)$.
Then there exists an element $u$ of the subgroup generated by $x_1,\dots,x_{n-1}$
such that $u(y_1\dots y_n)u^{-1} = x_1\dots x_n$.
\endproclaim

\demo{ Proof } Induction on $n$. If $n=1$, the statement is obvious.
Suppose that $n\ge 2$. Let $x=x_1\dots x_n$, $y=y_1\dots y_n$ and let
$G_k$ be the subgroup generated by $x_1,\dots, x_k$.
We have $y = A x_{n-1} B x_n C$ or $y = A x_n B x_{n-1} C$ with $A,B,C\in G_{n-2}$.
In the former case we have $y =  A x_{n-1} x_n BC$. In the latter case,
for $u_1 = x_{n-1}C\in G_{n-1}$, we have
$u_1 y u_1^{-1} = x_{n-1} CA x_n B = x_{n-1} x_n CAB$.
Thus, in the both cases, there exists $u_1\in G_{n-1}$ such that
$u_1 y u_1^{-1} = z_1\dots z_{n-1}$ where $(z_1,\dots, z_{n-1})$
is a permutation of $(x_1,\dots, x_{n-2}, x'_{n-1})$ for $x'_{n-1} = x_{n-1}x_n$.
Since $x'_{n-1}$ commutes with $x_1,\dots, x_{n-3}$, by the induction hypothesis
there exists $u_2\in G_{n-2}$ such that
$u_2 (z_1\dots z_{n-1}) u_2^{-1} = x_1\dots x_{n-2} x'_{n-1} = x$
and hence $u y u^{-1} = x$ for $u=u_2u_1$.
\qed\enddemo

\demo{ Proof of Theorem \thPermut }
(a).
First, note that $\Tproj(d,d-2)$ can be represented by a braid which satisfies these
conditions. For example, one can choose the braid (\eqT).

Thus, it is enough to show that if $Y$ is another braid satisfying the same conditions, then
the projective closure of $X$ and $Y$ coincide as knots in $\RP^3$.
To this end we show that the braids $X'=\Delta X^{-1}$ and $Y'=\Delta Y^{-1}$ are conjugate in this case.
Indeed, $X'$ and $Y'$ are permutation braids by [\refEM; Theorem 2.6] and
$$
    e(X')=e(Y')=e(\Delta)-e(X)={d(d-1)\over 2} - {(d-1)(d-2)\over 2} = d-1.
$$
Thus each of $X'$ and $Y'$ is a product of $d-1$ generators.
The connectedness of the projective braid closure implies that
each $\sigma_i$, $1\le i<d$, occurs in this product and we conclude
that $X'$ and $Y'$ are conjugate by Lemma \lemCoxeter.

\smallskip
(b). Proof is similar but simpler (we do not use [4; Theorem 2.6] and we apply Lemma \lemCoxeter\
directly to $A$).
\qed\enddemo



\head\sectExist. Construction of $MW_\lambda$-links of arbitrary genus (proof of Theorem \thExist)
\endhead

Theorem \thExist\ is an immediate consequence of Lemma~\lemExistOne\ and Lemma~\lemExistTwo\ below.
We say that a smooth irreducible algebraic curve $C$ in $\CP^n$ is {\it special} if the
the plane section divisor $D$ on $C$ is special, i.e., $h^1(D)>0$ where, as usual, $h^i(D)$ is an
abbreviation for $\dim H^i(C,\Cal O_C(D))$.

\proclaim{ Lemma \lemExistPrtrb }
Let $A$ be a smooth irreducible real algebraic non-special curve of genus $g$ and degree $d$ 
in $\RP^3$
and let $\ell$ be a real line which crosses $A$ at a single point $p$ without tangency.
Then there exists a smooth irreducible real algebraic non-special curve $A'$ of degree $d+1$ and genus $g$
arbitrarily close to $A\cup\ell$ whose real locus is isotopic to the link obtained
from $\R A\cup\R\ell$ by any of the two local modifications at $p$ shown in Figure \figSmoo.

Moreover, the curve $A'$ can be chosen isomorphic to $A$ as an abstract real algebraic curve. 
\endproclaim

\midinsert
\centerline{ \epsfxsize=50mm\epsfbox{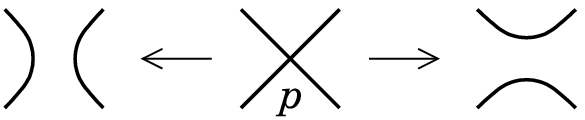}}
\botcaption{Figure \figSmoo }
Local smoothings in Lemma \lemExistPrtrb
\endcaption
\endinsert

\noindent{\bf Remark \remExist.}
\roster
\item
For $g=0$ this result follows from [\refB;  Theorem 2.4].
\smallskip\item
A simple example showing that Lemma \lemExistPrtrb\ may fail for special curves
can be found for $(d,g)=(4,3)$ (then $A$ is contained in a plane).
In this case $A\cup\ell$ cannot be smoothed out into a $5$-th degree curve of genus $3$ because
smooth spatial curves of degree 5 and genus 3 do not exist. Indeed, such a curve should sit either in
a plane or in an irreducible quadric by [\refMO; Proposition~1].
However the maximal genus of a quintic curve in a nonsingular quadric is
$\max_{a+b=5} (a-1)(b-1) = 2$ and the genus of any smooth quintic curve in a plane (resp. in a quadratic cone)
is $6$ (resp. $2$).
\smallskip\item
Lemma \lemExistPrtrb\ also admits an easy generalization to the case of two
non-special irreducible real algebraic curves $A$ and $B$ intersecting at a single
point $p$ without tangency. In such case there exists a smooth curve $C$ of
degree equal to the sum of degrees, and of genus equal to the sum of genera
such that its real locus is isotopic to the link obtained from $\R A\cup\R B$ by
any of the two local modifications at $p$ shown in Figure 8. However for the purposes
of this paper the special case given by Lemma \lemExistPrtrb\ suffices.
\endroster
\smallskip

\demo{ Proof }
Let us choose homogeneous coordinates $(x:y:z:w)$ so that $p$ is not on the coordinate planes and
$\ell$ passes through the point $(0:0:0:1)$.
Up to rescaling, the embedding of $A$ into $\P^3$ is determined by the effective divisors $D_x$, $D_y$,
$D_z$, $D_w$ which are cut on $A$ by the coordinate planes.
We are going to define $A'$ as the image of an embedding $f:A\to\P^3$ defined by new effective divisors
$D'_x$, $D'_y$, $D'_z$, $D'_w$. To achieve the result, it is enough to choose these divisors
so that they all belong to the same linear system and
$$
    D'_x = D_x+p,\quad D'_y=D_y+p,\quad D'_z=D_z+p,\quad D'_w = \wt D_w + \tilde p
$$
where $\wt D_w$ is close to $D_w$ and
$\tilde p$ is a real point on $A$ close
to $p$ and chosen on a required side of $p$ (the choice between the two local smoothings
in Figure \figSmoo\ is determined by the side from which $\tilde p$ approaches to $p$;
see Figure \figSmooG).

\midinsert
\centerline{ \epsfxsize=85mm\epsfbox{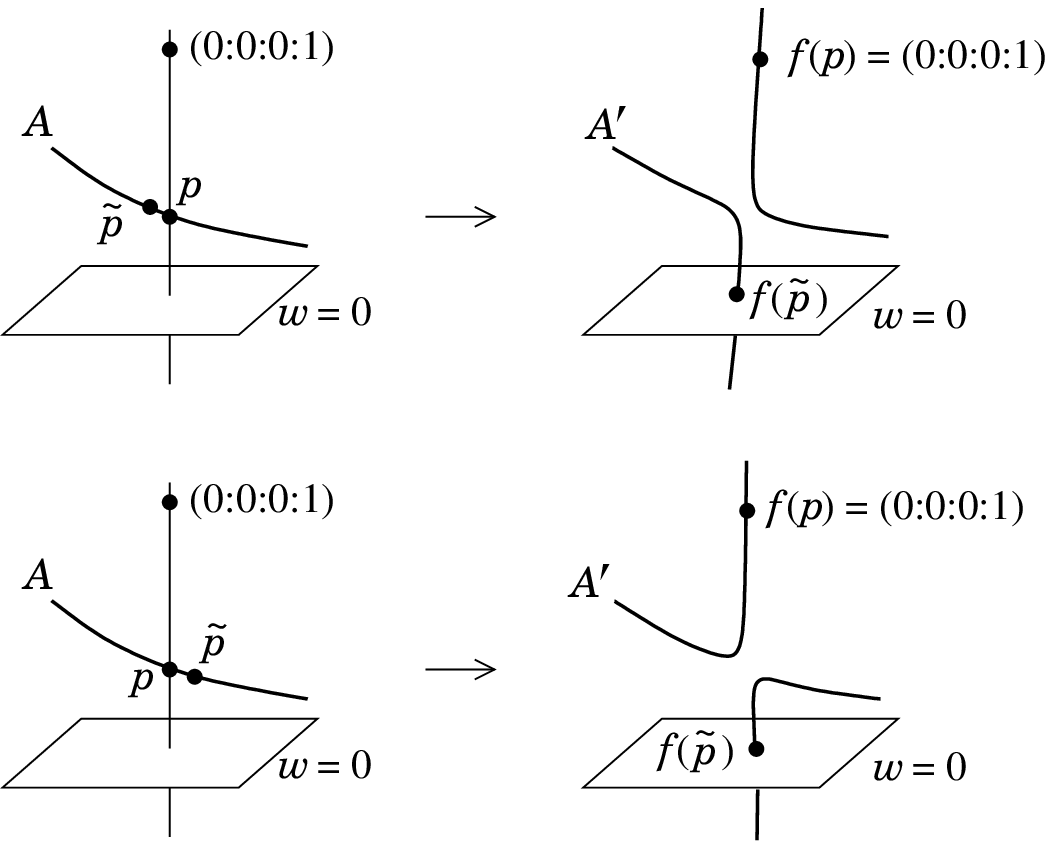}}
\botcaption{Figure \figSmooG }
Global smoothings in Lemma \lemExistPrtrb
\endcaption
\endinsert

For a positive integer $n$, let $s_n:S^n(A)\to\Cal J_A$ be the Abel-Jacobi mapping
where $S^n(A)$ is the $n$-th symmetric power of $A$
and $\Cal J_A$ is the Jacobean of $A$ (to define $s_n$ we need to fix a point $p_0$ on $A$).
Note that the both varieties admit natural real structure (an antiholomorphic involution Conj) and $s_n$
is a Conj-equivariant regular morphism (we assume here that $p_0$ is real).
The real locus $\R\Cal J_A$ is a subgroup of $\Cal J_A$. It is a union of
$2^{l-1}$ real $g$-dimensional tori where $l=b_0(\R A)$ (see [\refGroH], [\refN; \S8]).

It is clear that $D'_x\sim D'_y\sim D'_z$. The condition $D'_x\sim D'_w$ reads
$s_{d+1}(D_x+p) = s_{d+1}(\widetilde D_w+\tilde p)$. Since $D_x\sim D_w$, this is equivalent to
$$
      s_d(D_w) - s_d(\widetilde D_w) = s_1(\tilde p)-s_1(p).
$$
Thus we need to prove that for any sufficiently small $h\in\R\Cal J_A$ there exists
$\widetilde D_w\in\R S_d(A)$ close to $D_w$ such that
$s_d(D_w) - s_d(\widetilde D_w) = h$.
This follows from the surjectivity of the differential d$s_d$ at non-special
divisors $D$.

It remains to check that $A'$ is non-special.
Indeed, the corresponding plane section divisor on $A$ is $D'_x=D_x+p$. We have $h^0(K_A-D_x)=0$
because $A$ is non-special. Hence $h^0(K_A-D'_x)=0$.
\qed\enddemo

\proclaim{ Lemma \lemNonSpe }
Let $X=\P^1\times\P^1$
and let $C$ be a nonsingular algebraic curve on $X$ of bidegree $(2,g+1)$, $g\ge 0$.
Let $D$ be a divisor on $C$ which is cut by a generic
curve of bidegree $(1,1)$.
Then $h^1(D)=0$.
\endproclaim

\demo{ Proof }
Let $A$ and $B$ be algebraic curves on $X$ of bidegree $(1,0)$ and $(0,1)$ respectively, and let
Let $A|_C$ and $B|_C$ be the divisors in $C$ which are cut by them. 
By Serre's duality we have $h^1(D) = h^0(K_C - D)$. Let us show that $|K_C-D|$ is empty. Suppose it is not.
Let $D_1\in|K_C-D|$, thus
$$
      K_C\sim D_1+D = D_1 + A|_C + B|_C.                   \eqno(\eqNonSpe)
$$
The curve $C$ is hyperelliptic of genus $g$, and the hyperelliptic projection $\pi:C\to\P^1$ is the restriction
of the projection of $\P^1\times\P^1$ onto the first factor (this projection contracts $B$ to a single point).
Without loss of generality we may assume that $A$ is transverse to $C$, thus $\pi(\supp(A|_C))$ consists
of $A\cdot C=A\cdot(2A+(g+1)B)=g+1$ distinct points. Since all summands in the right hand side of (\eqNonSpe)
are effective divisors, we conclude that the support of $D_1+D$ contains that of $A|_C$ whence
$\pi(\supp(D_1+D))$ has at least $g+1$ points.

On the other hand, $D_1+D\in |K_C|$ (see (\eqNonSpe)) and it is well known (see, e.g.,
[\refGH; Ch.~3, \S 3.5]) that
any element of $|K_C|$ is supported by a union of at most $g-1$ fibers of $\pi$. The obtained contradiction
completes the proof.
 \qed\enddemo

\proclaim{ Lemma \lemExistOne }
For any $g\ge 0$ there exists an irreducible non-special real algebraic link
of genus $g$ and degree $d=g+3$, isotopic to $W_g(1,\dots,1)$.
\endproclaim

\demo{ Proof } Let $L_0,\dots,L_g$ be real lines belonging to one ruling
of a hyperboloid $H$, and let $L$ and $\bar L$ be a pair of complex conjugate
lines from the other ruling.
We smooth out (remaining on $H$) the union of all these $d=g+3$ lines.
The genus of the resulting curve is $g$. It is non-special by Lemma \lemNonSpe, and
its real locus is isotopic to $\R L_0\cup\dots\cup\R L_g$ which is a projective Hopf link.
\qed\enddemo

\proclaim{ Lemma \lemExistTwo }
Let $a_0,\dots,a_g$ be positive integers and $d=a_0+\dots+a_g+2$.
Suppose that $W_g(a_0,\dots,a_g)$ is realizable by an irreducible non-special
real algebraic link of degree $d$ and genus $g$.
Then $W_g(a_0+1,a_1,\dots,a_g)$
is realizable by an irreducible non-special real algebraic link of degree $d+1$ and genus $g$. \qed
\endproclaim

\demo{ Proof }
Let $\R A$ be a non-special real algebraic link of type $W_g(a_0,\dots,a_g)$.
We shall use for it the notation introduced in \S\sectHypSpa\ and \S\sectT.
Let $TK_0$ be the component of the tangent surface of $A$
corresponding to the $a_0$-component $K_0$ of $\R A$.
For a point $p$ close to $TK_0$ the curve $C_p\subset\pp^2$
(obtained as the image of $A$ under the projection $\pi_p:\pp^3\setminus\{p\}\to\pp^2$)
is a hyperbolic curve.

Since $p$ is close to $TK_0$, the curve $C_p$ is obtained by a perturbation
of a cuspidal hyperbolic curve $C_{p'}$, $p'\in TK_0$, with a cusp on
the projection of $K_0$. Thus we may choose $q\in\hyp(C_p)$ such that
$\ell=\pi^{-1}_p(q)\cap K_0\neq\varnothing$ so that $A$ and $\ell$ intersect
at a single point without tangency.

The two local modifications of $\R A\cup\R\ell$ provided by Lemma 5.1
produce non-special real algebraic links with the invariant $w_\lambda$
equal to $w_\lambda(\R A)\pm (d-1)$. Denote by $\R B$ the one corresponding
to the plus sign so that $w_\lambda(\R B)=N_d-g+(d-1)=N_{d+1}$.

By Theorem~\thAnyGenus\ the isotopy type of $\R B$ is $W_g(b_0,\dots,b_g)$ where
the numbers $b_j$ are determined by Lemma~\lemLK. We get $b_j=a_j$
for $j>0$ and $b_0=a_0+1$ once we keep the numeration of the components
of $\R B$ and $\R A$ consistent.
\qed\enddemo



\subhead Acknowledgement
\endsubhead
We are grateful to Oleg Viro for very useful and fruitful discussions.

\enddocument